\documentclass[10pt]{amsart}
\usepackage{amsmath,mathrsfs,cite}
\usepackage[colorlinks=true,urlcolor=blue,
citecolor=red,linkcolor=blue,linktocpage,pdfpagelabels,bookmarksnumbered,bookmarksopen]{hyperref}
\usepackage{color,mathrsfs}
\usepackage[english]{babel}
\usepackage[left=2.65cm,right=2.65cm,top=2.5cm,bottom=2.5cm]{geometry}

\numberwithin{equation}{section}
\newtheorem{theorem}{Theorem}[section]

\newtheorem{remark}[theorem]{Remark}

\theoremstyle{definition}

\newcommand{\R}{{\mathbb R}}

\newcommand{\dvg}{{\rm div}}

\title[Boundary behavior for a singular PDE]{Boundary behavior for a singular \\ quasi-linear elliptic equation}

\author[M.~Squassina]{Marco Squassina}

\address{Universit\`a degli Studi di Verona
\newline\indent
Strada Le Grazie 15, I-37134 Verona, Italy}
\email{marco.squassina@univr.it}

\thanks{Work partially supported by PRIN 
 {\em ``Variational and Topological Methods in the Study of Nonlinear Phenomena''}}


\begin{document}

\begin{abstract}
In a smooth bounded domain we obtain existence and uniqueness, regularity and 
boundary behavior for a class of singular quasi-linear elliptic equations.
\end{abstract}

\maketitle

\section{Introduction and result}

Let $\Omega\subset\R^N$ be a smooth bounded domain.
The aim of this note is to establish existence, uniqueness and 
boundary behavior of the solutions to the singular quasi-linear problem
\begin{equation}
\label{prob}
\begin{cases}
\,-\dvg(a(u)Du)+\frac{a'(u)}{2}|Du|^2=f(u)   & \text{in $\Omega$,} \\
\,\, u>0 &  \text{in $\Omega$,}  \\
\,\, u=0 &  \text{on $\partial\Omega$,}    
\end{cases}
\end{equation}
where $a:(0,+\infty)\to (0,+\infty)$ is a $C^1$ function 
bounded away from zero, $f:(0,+\infty)\to (0,+\infty)$ is a $C^1$ function
and there exist $a_0>0$ and $f_0>0$ such that
\begin{align}
\lim_{s\to 0^+} a(s)s^{2\mu}=a_0, & \qquad\text{for some $0\leq \mu<1$},   \label{ass-a} \\
\lim_{s\to 0^+} f(s)s^\gamma=f_0, & \qquad\text{for some $\gamma>1$},   \label{ass-f} \\
2f'(s)a(s)\leq f(s)a'(s), & \qquad\text{for every $s>0$}. \label{ass-fg}
\end{align}
We refer to \cite{squassmono} and to the references included
for the interest and motivations to analyze these equations.
In the one dimensional case, equations such as \eqref{prob} typically arise in certain problems in 
fluid mechanics and pseudo-plastic flow (see e.g.\ \cite{callegari,stuart}). In the semi-linear case $a\equiv 1$
various results about existence, uniqueness and asymptotic behavior of the solutions have been obtained
in the literature so far (see \cite{craRabTar,lairshak,lazmec,zhangcheng,gherad-paper}, the monographs \cite{gherad-book,radu-hand}
and the references therein). 
Under assumptions \eqref{ass-a}-\eqref{ass-fg}, if ${\rm d}(x,\partial\Omega)$ denotes the distance of a point $x$ in $\Omega$
from the boundary $\partial\Omega$, we shall prove the following

\begin{theorem}
	\label{main}
Problem~\eqref{prob} has a unique solution $u\in C^2(\Omega)\cap C(\bar\Omega)$ and
there exist $\Gamma,\Gamma',\Gamma''>0$ with
\begin{align*}
\Gamma{\rm d}(x,\partial\Omega)^{{\frac{2}{1+\gamma-2\mu}}}\leq u(x)\leq \Gamma'
{\rm d}(x,\partial\Omega)^{{\frac{2}{1+\gamma-2\mu}}}, \quad
|Du(x)|\leq \Gamma''{\rm d}(x,\partial\Omega)^{\frac{1-\gamma+2\mu}{1+\gamma-2\mu}},&\quad\,\, \text{as 
${\rm d}(x,\partial\Omega)\to 0$.}  
\end{align*}
Moreover $u\in {\rm Lip}(\bar\Omega)$ if $1<\gamma\leq 1+2\mu$,
$u\in C^{0,\frac{2}{1+\gamma-2\mu}}(\bar\Omega)$ if $\gamma>1+2\mu$
and $u\in H^1_0(\Omega)$ if $\gamma<3-2\mu$.
\end{theorem}

\noindent
Whence, in some sense, the functions $a$ and $f$ compete for the vanishing rate of the solution and for its gradient upper bound,
near the boundary $\partial\Omega$. Furthermore, the range for Lipschitz continuity of $u$ up to the boundary is 
$\gamma\leq 1+2\mu$, $\gamma\neq 1$ thus enlarged with respect to the one for the semi-linear case, namely $0<\gamma<1$,
see also Remark~\ref{globalrem}.
In \cite{GlasquII}, the author jointly with F.\ Gladiali have recently performed a complete study about existence
and qualitative behavior around $\partial\Omega$ of the solutions to the problem
\begin{equation*}
\begin{cases}
\,\dvg(a(u)Du)-\frac{a'(u)}{2}|Du|^2=f(u)   & \text{in $\Omega$,} \\
\noalign{\vskip2pt}
\,\text{$u(x)\to+\infty$\quad as ${\rm d}(x,\partial\Omega)\to 0$,}    &
\end{cases}
\end{equation*}
covering situations where $a$ and $f$ have an exponential, polynomial or logarithmic 
type growth at infinity and a nonsingular behavior around the origin. On the contrary,
here we focus on the singular behavior at
the origin for $a$ and $f$ with the action of the source $f$ being in some sense predominant at zero upon the
diffusion $a$, due to the constraint $\gamma>1>\mu$. 
Without loss of generality we assume that $a$ grows as $s^k$ and $f$ decays as $s^{-p}$ as $s\to+\infty$
for some $k\geq 0$ and $p>1$, in which case one can also obtain estimates for $u$ and $|Du|$ valid
on the whole $\bar\Omega$, as pointed out in  Remark~\ref{globalrem}.
In the particular case when $a\equiv 1$, problem \eqref{prob} reduces to $-\Delta u =f(u)$
and the above estimates reduce to $\Gamma{\rm d}(x,\partial\Omega)^{2/(1+\gamma)}\leq u(x)\leq \Gamma'
{\rm d}(x,\partial\Omega)^{2/(1+\gamma)}$ and $|Du(x)|\leq \Gamma''{\rm d}(x,\partial\Omega)^{(1-\gamma)/(1+\gamma)}$
as ${\rm d}(x,\partial\Omega)\to 0$, consistently with the results of \cite{craRabTar,zhangcheng}.
Following the line of \cite{lairshak,zhangcheng}, some easy adaptations 
of Theorem~\ref{main} can be obtained to cover the case of non-autonomous nonlinearities
such as $q(x)f(u)$ in place of $f(u)$ and of unbounded domains of $\R^N$. We leave these further 
developments to the interested reader.
As an example of $f$ and $a$ satisfying \eqref{ass-a}-\eqref{ass-fg} one can take
$f(s)=s^{-\gamma}$, $a(s)=s^{-2\mu}$ for $s\leq s_0$ and $a(s)=\theta(s)$ for $s\geq s_0$
for some $s_0>0$, with $\theta\in C^1\cap L^\infty$ bounded away from zero, 
$\theta(s_0)=s_0^{-2\mu}$, $\theta'(s_0)=-2\mu s_0^{-2\mu-1}$
and $s\theta'(s)+2\gamma \theta(s)\geq 0$ for $s\geq s_0$.

\section{Proof of the result}
\noindent

In this section, we prove Theorem~\ref{main}. We shall assume that conditions \eqref{ass-a}-\eqref{ass-fg} hold.
In order to get information about existence, uniqueness and the boundary behavior of the solutions to \eqref{prob}, 
we convert the quasi-linear problem \eqref{prob} into a corresponding semi-linear
problem through a change of variable procedure involving the Cauchy problem for $g\in C^2((0,+\infty))\cap C([0,+\infty))$,
\begin{equation}
	\label{cauchy}
	\begin{cases}
g'(s)=\frac{1}{\sqrt{a(g(s))}},\quad \text{for $s>0$}, &\\
g(0)=0, & \\
\noalign{\vskip3pt}
g(s)>0, \quad\text{for $s>0$}. &
\end{cases}
\end{equation}
Due to the requirement $g>0$, the solutions of~\eqref{cauchy} are unique and solve $\int_0^{g(s)}\sqrt{a(\xi)}d\xi=s$, for $s>0$.
The solution is global for $s>0$ since $a$ is bounded away from zero.
This procedure was also followed in  \cite{GlasquII} in the framework of explosive solutions, although
there $g$ is $C^2$ around the origin and defined on $\R$. Now, since $g\in C^2((0,+\infty))\cap C([0,+\infty))$
and it is strictly increasing, it is readily seen by a direct computation that 
$u\in C^2(\Omega)\cap C(\bar\Omega)$ is a positive solution to \eqref{prob} if and only if 
$v=g^{-1}(u)\in C^2(\Omega)\cap C(\bar\Omega)$ is a positive solution to 
$-\Delta v=h(v)$ in $\Omega$, where we have set $h(s):=f(g(s))/\sqrt{a(g(s))}$ for $s>0$.
Let us now obtain the asymptotic behavior of the solution $g$ to problem~\eqref{cauchy} as $s\to 0^+$ depending 
of the assigned asymptotic behavior of $a$ as $s\to 0^+$, given by~\eqref{ass-a}. For every $0\leq \mu<1$, we have 
\begin{equation}
	\label{gatzero}
\lim_{s\to 0^+}\frac{g(s)}{s^{\frac{1}{1-\mu}}}=(1-\mu)^{\frac{1}{1-\mu}} a_0^{\frac{1}{2(\mu-1)}}.      
\end{equation}
In fact, taking into account \eqref{cauchy}, by l'H\^ opital's rule we have
\begin{align*}
\lim_{s\to 0^+}\frac{g(s)}{s^{\frac{1}{1-\mu}}}&=(1-\mu)\lim_{s\to 0^+}\frac{g'(s)}{s^{\frac{\mu}{1-\mu}}}
=(1-\mu)\lim_{s\to 0^+}\frac{1}{s^{\frac{\mu}{1-\mu}}\sqrt{a(g(s))}} \\
& =(1-\mu)\lim_{s\to 0^+}\frac{g^\mu(s)}{s^{\frac{\mu}{1-\mu}}\sqrt{a(g(s))g^{2\mu}(s)}} =
\frac{(1-\mu)}{\sqrt{a_0}} \Big(\lim_{s\to 0^+}  \frac{g(s)}{s^{\frac{1}{1-\mu}}}\Big)^\mu,
\end{align*}
which yields the claim. Moreover, by virtue of~\eqref{ass-a}, \eqref{ass-f} and \eqref{gatzero}, we have
\begin{equation}
	\label{accat0}
\lim_{s\to 0^+}\frac{h(s)}{s^{\frac{\mu-\gamma}{1-\mu}}}=
\lim_{s\to 0^+} f(g(s))g(s)^{\gamma}\frac{1}{\sqrt{a(g(s))g(s)^{2\mu}}}
\Big[\frac{g(s)}{s^{\frac{1}{1-\mu}}}\Big]^{\mu-\gamma}=
f_0 a_0^{\frac{1-\gamma}{2(\mu-1)}}(1-\mu)^{-\frac{\gamma-\mu}{1-\mu}}.
\end{equation}
Observe also that, since $a(s)\sim a_\infty s^k$ as $s\to+\infty$
and $f(s)\sim f_\infty s^{-p}$ as $s\to+\infty$ for some $a_\infty,f_\infty>0$, $k\geq 0$ and $p>1$,
if $g$ still denotes the solution to \eqref{cauchy}, we have three facts:
	\begin{equation}
		\label{limzze}
		\int_1^{+\infty}\frac{f(g(s))}{\sqrt{a(g(s))}}ds<+\infty,  \quad
\lim_{s\to 0^+}\frac{f(g(s))}{\sqrt{a(g(s))}}=+\infty,\quad
\text{$s\mapsto \frac{f(g(s))}{\sqrt{a(g(s))}}$ is nonincreasing.}
	\end{equation} 
The first property follows immediately from the limit
\begin{equation*}
\lim_{s\to +\infty}\frac{g(s)}{s^{\frac 2{k+2}}}=\big(\frac{k+2} 2 \frac 1{\sqrt a_\infty} \big)^{\frac 2{k+2}},
\end{equation*}
which was proved in \cite{GlasquII}. The other properties follow by \eqref{accat0} and
\eqref{ass-fg}, respectively. By virtue of \eqref{gatzero} we now prove that, 
for every $\gamma>1$ and $0\leq \mu<1$, there holds
\begin{equation}
	\label{limfracintt}
\lim_{s\to 0^+}\frac{\int_{g(s)}^{+\infty} f(\xi)d\xi}{s^{\frac{1-\gamma}{1-\mu}}}=f_0 a_0^{\frac{\gamma-1}{2(1-\mu)}}(\gamma-1)^{-1}
(1-\mu)^{-\frac{\gamma-1}{1-\mu}}.
\end{equation}
In fact, since $\int_{g(s)}^{+\infty} f(\xi)d\xi\to +\infty$ for all $\gamma>1$, it follows
\begin{align*}
\lim_{s\to 0^+}\frac{\int_{g(s)}^{+\infty} f(\xi)d\xi}{s^{\frac{1-\gamma}{1-\mu}}}&=
\frac{\mu-1}{1-\gamma}\lim_{s\to 0^+}\frac{f(g(s))}{s^{\frac{\mu-\gamma}{1-\mu}}\sqrt{a(g(s))}}
=\frac{\mu-1}{1-\gamma}\lim_{s\to 0^+}\frac{f(g(s))g^\mu(s)}{s^{\frac{\mu-\gamma}{1-\mu}}\sqrt{a(g(s))g^{2\mu}(s)}} \\
& =\frac{\mu-1}{(1-\gamma)\sqrt{a_0}}\lim_{s\to 0^+}f(g(s))g^\gamma(s)
\lim_{s\to 0^+}\frac{g^{\mu-\gamma}(s)}{s^{\frac{\mu-\gamma}{1-\mu}}} \\
& =\frac{\mu-1}{(1-\gamma)}\frac{f_0}{\sqrt{a_0}}
\lim_{s\to 0^+}\frac{g^{\mu-\gamma}(s)}{s^{\frac{\mu-\gamma}{1-\mu}}} 
=\frac{\mu-1}{(1-\gamma)}\frac{f_0}{\sqrt{a_0}}
\Big(\lim_{s\to 0^+}\frac{g(s)}{s^{\frac{1}{1-\mu}}}\Big)^{\mu-\gamma} \\
& =f_0 a_0^{\frac{\gamma-1}{2(1-\mu)}}(\gamma-1)^{-1}
(1-\mu)^{-\frac{\gamma-1}{1-\mu}}.
\end{align*}
Now, for every $\ell\geq 0$, there exists a 
unique solution $\phi\in C([0,+\infty))\cap C^2((0,+\infty))$ of the problem
\begin{equation}
	\begin{cases}
	\label{asympt-prob}
	\displaystyle
\phi'(s)=\sqrt{\ell^2+2\int_{g(\phi(s))}^{+\infty}f(\xi)d\xi}, \qquad \text{for $s>0$}, & \\
\noalign{\vskip4pt}
\, \phi(0)=0, \quad \phi(s)>0, \qquad\text{for $s>0$}, & \\
\noalign{\vskip6pt}
\lim\limits_{s\to+\infty}\phi'(s)=\ell,\,\,\,\,  \lim\limits_{s\to+\infty}\phi(s)=+\infty. & 
\end{cases}
\end{equation}
To prove this, taking into account \eqref{limzze}, it is sufficient to apply \cite[Lemma 1.3]{zhangcheng}.
Notice that, in particular, the solutions to problem \eqref{asympt-prob} locally (namely for every fixed $a>0$) 
solve the second order problem
\begin{equation}
	\begin{cases}
	\label{asympt-prob-ii}
	\displaystyle
-\phi''(s)=\frac{f(g(\phi(s))}{\sqrt{a(g(\phi(s)))}}, \qquad & \text{for $0<s\leq a$}, \\
\noalign{\vskip4pt}
\, \phi(0)=0, \quad \phi(s)>0, & \text{for $0<s\leq a$}.
\end{cases}
\end{equation}
We can now prove that, for every $\gamma>1$, $0\leq\mu<1$ and $\ell\geq 0$, there holds
\begin{equation}
	\label{zerobeh}
\lim_{s\to 0^+}\frac{\phi(s)}{s^{\frac{2-2\mu}{1+\gamma-2\mu}}}=
\Big(\frac{1+\gamma-2\mu}{2-2\mu}\Big)^{\frac{2-2\mu}{1+\gamma-2\mu}}		
f_0^{\frac{1-\mu}{1+\gamma-2\mu}} a_0^{\frac{\gamma-1}{2(1+\gamma-2\mu)}}(\gamma-1)^{\frac{\mu-1}{1+\gamma-2\mu}} 
	(1-\mu)^{\frac{1-\gamma}{1+\gamma-2\mu}},
\end{equation}
where $\phi$ denotes the unique solution to \eqref{asympt-prob}. 
In fact, by l'H\^ opital's rule and \eqref{limfracintt}, we obtain
	\begin{align*}
	\lim_{s\to 0^+}\frac{\phi(s)}{s^{\frac{2-2\mu}{1+\gamma-2\mu}}} &=
	\lim_{s\to 0^+}\Big[\frac{\phi(s)^{\frac{1+\gamma-2\mu}{2-2\mu}}}{s}\Big]^{\frac{2-2\mu}{1+\gamma-2\mu}} 
	=\lim_{s\to 0^+}\Bigg[\frac{\phi(s)^{\frac{1+\gamma-2\mu}{2-2\mu}}}{
	\displaystyle\int_0^{\phi(s)}\frac{d\tau}{\sqrt{\ell^2+\int_{g(\tau)}^{+\infty}f(\xi)d\xi}}  }\Bigg]^{\frac{2-2\mu}{1+\gamma-2\mu}} \\
&=\Big(\frac{1+\gamma-2\mu}{2-2\mu}\Big)^{\frac{2-2\mu}{1+\gamma-2\mu}}
	\lim_{s\to 0^+}\Big[\phi(s)^{\frac{\gamma-1}{2-2\mu}}\sqrt{\ell^2+\int_{g(\phi(s))}^{+\infty}f(\xi)d\xi}\Big]^{\frac{2-2\mu}{1+\gamma-2\mu}}    \\
	&=\Big(\frac{1+\gamma-2\mu}{2-2\mu}\Big)^{\frac{2-2\mu}{1+\gamma-2\mu}}
		\lim_{s\to 0^+}\left(\frac{\ell^2+\int_{g(\phi(s))}^{+\infty}f(\xi)d\xi}{ \phi(s)^{\frac{1-\gamma}{1-\mu}} }\right)^{\frac{1-\mu}{1+\gamma-2\mu}}   \\
		&=\Big(\frac{1+\gamma-2\mu}{2-2\mu}\Big)^{\frac{2-2\mu}{1+\gamma-2\mu}}
			\lim_{s\to 0^+}\left(\frac{\int_{g(s)}^{+\infty}f(\xi)d\xi}{s^{\frac{1-\gamma}{1-\mu}} }\right)^{\frac{1-\mu}{1+\gamma-2\mu}}  \\
			\noalign{\vskip4pt}
		&=\Big(\frac{1+\gamma-2\mu}{2-2\mu}\Big)^{\frac{2-2\mu}{1+\gamma-2\mu}}		
		f_0^{\frac{1-\mu}{1+\gamma-2\mu}} a_0^{\frac{\gamma-1}{2(1+\gamma-2\mu)}}(\gamma-1)^{\frac{\mu-1}{1+\gamma-2\mu}} 
			(1-\mu)^{\frac{1-\gamma}{1+\gamma-2\mu}}.
	\end{align*}
We are now ready to conclude the proof 
of Theorem~\ref{main}. 
In light of \cite[Theorem 1.1]{craRabTar}, since $h(s)\to +\infty$ as $s\to 0^+$ by \eqref{limzze}
and $h$ is non-increasing for $s>0$ by \eqref{ass-fg}, there exists a unique positive solution $z\in C^2(\Omega)\cap C(\bar\Omega)$
to $-\Delta v=h(v)$. Then $g(z)\in C^2(\Omega)\cap C(\bar\Omega)$ is a positive solution to \eqref{prob}. Assume that
$u_1,u_2\in C^2(\Omega)\cap C(\bar\Omega)$, $u_1,u_2>0$ solve \eqref{prob}. Then $g^{-1}(u_1),g^{-1}(u_2)>0$ solve
$-\Delta v=h(v)$. By uniqueness, we deduce $g^{-1}(u_1)=g^{-1}(u_2)$, in turn yielding $u_1=u_2$. By virtue of 
\cite[Theorem 2.2 and Theorem 2.5]{craRabTar} for any solution $v$ of $-\Delta v=h(v)$
there exist four constants $\Lambda_1,\Lambda_2,\Lambda_3,\Lambda_4>0$ such that
for ${\rm d}(x,\partial\Omega)$ small enough, $\Lambda_1 \phi({\rm d}(x,\partial\Omega))\leq 
v(x)\leq \Lambda_2 \phi({\rm d}(x,\partial\Omega))$ and
$|Dv(x)|\leq \Lambda_3\big[{\rm d}(x,\partial\Omega)h(\Lambda_4\phi({\rm d}(x,\partial\Omega)))+\phi({\rm d}(x,\partial\Omega))/{\rm d}(x,\partial\Omega)\big]$, 
being $\phi$ a solution to \eqref{asympt-prob-ii}. On account of formula~\eqref{zerobeh}, we can find two constants $\Theta_1,\Theta_2>0$
such that, for ${\rm d}(x,\partial\Omega)$ small enough
\begin{equation}
	\label{phicontrol}
\Theta_1 {\rm d}(x,\partial\Omega)^{{\frac{2-2\mu}{1+\gamma-2\mu}}}\leq \phi({\rm d}(x,\partial\Omega))\leq
\Theta_2 {\rm d}(x,\partial\Omega)^{{\frac{2-2\mu}{1+\gamma-2\mu}}},
\end{equation}
which yield in turn
$$
\Lambda_1\Theta_1 {\rm d}(x,\partial\Omega)^{{\frac{2-2\mu}{1+\gamma-2\mu}}}\leq 
v(x)\leq \Lambda_2\Theta_2 {\rm d}(x,\partial\Omega)^{{\frac{2-2\mu}{1+\gamma-2\mu}}}.
$$
Finally, since $g$ is increasing and $u=g(v)$, we have
$$
g\big(\Lambda_1\Theta_1 {\rm d}(x,\partial\Omega)^{{\frac{2-2\mu}{1+\gamma-2\mu}}}\big)\leq 
u(x)\leq g\big(\Lambda_2\Theta_2 {\rm d}(x,\partial\Omega)^{{\frac{2-2\mu}{1+\gamma-2\mu}}}\big).
$$
Finally, using \eqref{gatzero}, we obtain the desired controls on $u$.
Now, from \eqref{phicontrol}, as ${\rm d}(x,\partial\Omega)$ is small,
$$
\phi({\rm d}(x,\partial\Omega))/{\rm d}(x,\partial\Omega)\leq \Theta_2 {\rm d}(x,\partial\Omega)^{-{\frac{\gamma-1}{1+\gamma-2\mu}}}.
$$
On account of \eqref{accat0} and \eqref{zerobeh}, there exists a constant $\Theta_3>0$ such that 
\begin{align*}
& {\rm d}(x,\partial\Omega)h(\Lambda_4\phi({\rm d}(x,\partial\Omega))) =\\
& =\Lambda_4^{(\mu-\gamma)/(1-\mu)}{\rm d}(x,\partial\Omega)
\Big[\frac{h(\Lambda_4\phi({\rm d}(x,\partial\Omega)))}{\Lambda_4^{(\mu-\gamma)/(1-\mu)}
\phi^{(\mu-\gamma)/(1-\mu)}({\rm d}(x,\partial\Omega))}\Big]\phi^{(\mu-\gamma)/(1-\mu)}({\rm d}(x,\partial\Omega)) \\
& \leq \Theta_3{\rm d}(x,\partial\Omega){\rm d}(x,\partial\Omega)^{{\frac{2\mu-2\gamma}{1+\gamma-2\mu}}}=
\Theta_3 {\rm d}(x,\partial\Omega)^{-{\frac{\gamma-1}{1+\gamma-2\mu}}},
\end{align*}
for ${\rm d}(x,\partial\Omega)$ small enough, yielding 
$|Dv(x)|\leq \Lambda_3 (\Theta_2+\Theta_3){\rm d}(x,\partial\Omega)^{-{\frac{\gamma-1}{1+\gamma-2\mu}}}$.
Then, we have
\begin{align*}
|Du(x)| &=g'(v(x))|Dv(x)|=\frac{|Dv(x)|}{\sqrt{a(g(v(x)))}}=\frac{|Dv(x)|}{\sqrt{a(g(v(x)))g^{2\mu}(v(x))}}g^{\mu}(v(x)) &\\
& =\frac{1}{\sqrt{a(g(v(x)))g^{2\mu}(v(x))}}\frac{g^{\mu}(v(x))}{v(x)^{\mu/(1-\mu)}}v(x)^{\mu/(1-\mu)}|Dv(x)|
\leq\omega_1 v(x)^{\mu/(1-\mu)}|Dv(x)| \\
\noalign{\vskip4pt}
&\leq \omega_2{\rm d}(x,\partial\Omega)^{{\frac{2\mu}{1+\gamma-2\mu}}}{\rm d}(x,\partial\Omega)^{-{\frac{\gamma-1}{1+\gamma-2\mu}}}
=\omega_2{\rm d}(x,\partial\Omega)^{{\frac{1-\gamma+2\mu}{1+\gamma-2\mu}}}.
\end{align*}
for some $\omega_1,\omega_2>0$. In particular, if $1<\gamma\leq 1+2\mu$, it follows that $u$
is Lipschitz continuous up the boundary. If instead $\gamma>1+2\mu$,
by the above estimates for $u$ and $|Du|$, we find $\Theta_4>0$ such that
\begin{align*}
|D u^{\frac{1+\gamma-2\mu}{2}}(x)|& =\frac{1+\gamma-2\mu}{2}u^\frac{\gamma-1-2\mu}{2}(x)|Du(x)| \\
& \leq \Theta_4 {\rm d}(x,\partial\Omega)^{{\frac{\gamma-1-2\mu}{1+\gamma-2\mu}}}{\rm d}(x,\partial\Omega)^{{\frac{1-\gamma+2\mu}{1+\gamma-2\mu}}}
=\Theta_4,
\end{align*}
whenever ${\rm d}(x,\partial\Omega)$ is small enough. In turn, since $u^{\frac{1+\gamma-2\mu}{2}}$
is Lipschitz continuous, $0<2/(1+\gamma-2\mu)<1$ and $u=(u^{(1+\gamma-2\mu)/2})^{2/(1+\gamma-2\mu)}$ it
follows that $u$ is H\"older continuous up to the boundary $\partial\Omega$ with exponent $2/(1+\gamma-2\mu)$,
as desired.
Finally, concerning the Sobolev regularity of the solution $u$, observe that in light of 
\cite[Theorem 1.3-J2]{zhangcheng}, a necessary and sufficient condition for
$v$ to belong to $H^1_0(\Omega)$ is that
$$
\lim_{s\to 0^+}\int_s^1 \phi(\xi)h(\phi(\xi))d\xi<+\infty,
$$
and this, since by \eqref{accat0} and \eqref{zerobeh} 
$ \phi(\tau)h(\phi(\tau))\sim \tau^{\frac{2-2\gamma}{1+\gamma-2\mu}}$ as $\tau\to 0^+$, is satisfied
if and only if $\gamma<3-2\mu$. In turn, $u=g(v)\in H^1_0(\Omega)$ if $\gamma<3-2\mu$ since
$g$ is Lipschitz continuous on $(0,+\infty)$ being $a$ bounded away from zero. This concludes the proof of the theorem. \qed

\begin{remark}\rm
	\label{lipremark}
We know from Theorem~\ref{main} that a solution $u$ exists unique and it is Lipschitz continuous up to $\partial\Omega$
provided that $1<\gamma\leq 1+2\mu$. On the other hand, if $\gamma<1$, combining \eqref{accat0} with \cite[Theorem 2.25]{craRabTar}
any solution $v$ of $-\Delta v=h(v)$ is Lipschitz continuous, and so is $u$, since $u=g(v)$ and $g$ is
Lipschitz continuous on $(0,+\infty)$ since $a$ is bounded away from zero.
\end{remark}

\begin{remark}\rm
\label{globalrem}
Using \cite[Theorem 1.3]{zhangcheng} in place of \cite[Theorem 2.2 and Theorem 2.5]{craRabTar} we could also
state some global estimates for $u$ and $|Du|$ which are valid on the whole $\bar\Omega$ and not only in a small neighborhood
of the boundary $\partial\Omega$. Precisely, under the assumptions of Theorem~\ref{main}, 
there exist $\Lambda_1,\Lambda_2,\Lambda_3,\Lambda_4>0$ with
\begin{align*}
 & g(\Lambda_1 \phi({\rm d}(x,\partial\Omega)))\leq 
u(x)\leq  g(\Lambda_2 \phi({\rm d}(x,\partial\Omega))), \,\,\,\,\,\,\text{for $x\in\bar\Omega,$}\\
& |Du(x)|\leq \Lambda_3\frac{\big[{\rm d}(x,\partial\Omega)h(\Lambda_4\phi({\rm d}
(x,\partial\Omega)))+\phi({\rm d}(x,\partial\Omega))/{\rm d}(x,\partial\Omega)\big]}{\sqrt{a(u(x))}},
\,\,\,\,\,\,\text{for $x\in\bar\Omega,$}
\end{align*}
where $\phi$ denotes the solution to problem~\eqref{asympt-prob}. These formulas
are obtained though the monotonicity of $g$ and from $|Du(x)|=|Dv(x)|/\sqrt{a(u(x))}$ for all $x\in\bar\Omega$,
following by the relation $u=g(v)$.
\end{remark}

\begin{remark}\rm 
	\label{convexity}
Let $v$ be a smooth positive solution to 
$-\Delta v=h(v)$ in $\Omega$, $v=0$ on $\partial\Omega$,
	where $h>0$ is as in the proof of Theorem~\ref{main}. 
	Let us consider the map $\psi:[0,+\infty)\to [0,+\infty)$ defined by
	$$
	\psi(s):=\int_0^s \frac{1}{\sqrt{2{\mathcal H}(\xi)}}d\xi,\qquad
	{\mathcal H}(s):=\int_s^{+\infty} h(\xi)d\xi,
	\,\,\quad\text{for $s>0$.}
	$$
	It follows that $\psi'(s)=1/{\sqrt{2{\mathcal H}(s)}}>0$,
	$\psi''(s)={h(s)}/{(2{\mathcal H}(s))^{3/2}}>0$
	and $\psi$ satisfies $\psi''(s)=(\psi'(s))^3 h(s)$ for all $s>0$.
	Put $w:=-\psi(v)<0$, we have $|Dv|^2=|Dw|^2/(\psi'(v))^2$ as well as
	$\Delta w =-\psi''(v)|Dv|^2+\psi'(v)h(v)=\psi'(v) h(v)(1-|Dw|^2)$, yielding
	\begin{equation*}
	\Delta w=f(w,Dw),\qquad  f(w,Dw):=-\psi'(\psi^{-1}(-w)) h(\psi^{-1}(-w))(1-|Dw|^2),
	\end{equation*}
	with $w=0$ on $\partial\Omega$. Whenever $f>0$ and $z\mapsto 1/f(z,\cdot)$ is convex, one typically
	obtains some convexity of $w$ if $\Omega$ is convex (see \cite{kawohl}) and in turn some 
	convexity of superlevels of $u$ since $u=g\circ\psi^{-1}(-w)$ and since 
	$s\mapsto g\circ\psi^{-1}(s)$ is strictly increasing. See \cite[Sec.\ 3]{frapor}
	for the particular case $a\equiv 1$ and $f(s)=s^{-\gamma}$.
	\end{remark}

\medskip


\bigskip


\begin{thebibliography}{99}

\bibitem{frapor}
{\sc S.\ Berhanu, F.\ Gladiali, G.\ Porru}, 
Qualitative properties of solutions to elliptic singular problems,
{\em J. Inequal. Appl.} {\bf 3} (1999), 313--330.

\bibitem{craRabTar}
{\sc M.G. Crandall, P.H.\  Rabinowitz, L.\ Tartar}, 
On a Dirichlet problem with a singular nonlinearity, 
{\em Comm. Partial Differential Equations} {\bf 2} (1977), 193--222.

\bibitem{gherad-book}
{\sc M.\ Ghergu, V.\ Radulescu}, 
Singular Elliptic Problems. Bifurcation and Asymptotic Analysis, 
Oxford Lecture Series in Mathematics and Its Applications, {\bf 37} 
Oxford University Press, 320 pages, 2008.

\bibitem{gherad-paper}
{\sc M.\ Ghergu, V.\ Radulescu}, Multiparameter bifurcation and asymptotics for the
singular Lane-Emden-Fowler equation with a convection term, 
{\em Proceedings of the Royal Society of Edinburgh A} {\bf 135} (2005), 61--84.

\bibitem{GlasquII}
{\sc F.\ Gladiali, M. Squassina}, 
On explosive solutions for a class of quasi-linear elliptic equations, {\em preprint}.





\bibitem{kawohl}
{\sc B.\ Kawohl},
Rearrangements and convexity of level sets in PDE,
{\em Lecture Notes in Mathematics} {\bf 1150} Springer, Berlin.


\bibitem{lairshak}
{\sc A.V.\ Lair, A.W.\ Shaker}, 
Classical and weak solutions of a singular semilinear elliptic problem, 
{\em J. Math. Anal. Appl.} {\bf 211} (1997), 371--385.

\bibitem{lazmec}
{\sc A.C.\ Lazer, P.J.\ McKenna}, 
On a singular nonlinear elliptic boundary-value problem, 
{\em Proc. Amer. Math. Soc.} {\bf 111} (1991), 721--730.

\bibitem{callegari}
{\sc A.\ Nachman, A.\ Callegari}, 
A nonlinear singular boundary value problem in the theory of pseudoplastic fluids, 
{\em SIAM J. Appl. Math.} {\bf 38} (1980), 275--281. 

\bibitem{radu-hand}
{\sc V.\ Radulescu},
Singular phenomena in nonlinear elliptic problems. From blow-up boundary 
solutions to equations with singular nonlinearities, in Handbook of 
Differential Equations: Stationary Partial Differential Equations, {\bf 4} 
(Michel Chipot, Editor), North-Holland Elsevier Science, Amsterdam, 2007, 483--591.

\bibitem{squassmono}
{\sc M. Squassina},  
Existence, multiplicity, perturbation, and concentration results for a class of quasi-linear elliptic problems,
{\em Electron. J. Differential Equations}, Monograph {\bf 7} (2006), 213 pages.

\bibitem{stuart}
{\sc C.A.\ Stuart}, 
Existence theorems for a class of nonlinear integral equations, 
{\em Math. Z.} {\bf 137} (1974), 49--66.

\bibitem{zhangcheng}
{\sc Zhijun Zhang,  Jiangang Cheng}, 
Existence and optimal estimates of solutions for singular nonlinear Dirichlet problems,
{\em Nonlinear Anal.} {\bf 57} (2004), 473--484.

\end{thebibliography}
\end{document}